\documentclass{amsart}
\usepackage{amssymb,latexsym}
\input epsf
\newtheorem{theorem}{Theorem}[section]

\newtheorem{corollary}[theorem]{Corollary}

\newtheorem{remark}[theorem]{Remark}

\newtheorem{example}[theorem]{Example}

\newcommand{\NN}{{\mathbb N}}
\newcommand{\ZZ}{{\mathbb Z}}
\newcommand{\QQ}{{\mathbb Q}}
\newcommand{\RR}{{\mathbb R}}
\newcommand{\KK}{{\mathbb K}}
\newcommand{\aA}{{\mathcal A}}

\newcommand{\fF}{{\mathcal F}}
\newcommand{\gG}{{\mathcal G}}
\newcommand{\hH}{{\mathcal H}}

\renewcommand{\to}{\rightarrow}

\newcommand{\im}{{\rm Im}}

\newcommand{\vol}{{\rm vol}}

\begin{document}
\title[A reciprocity theorem for arrangements]
{A combinatorial reciprocity theorem for hyperplane arrangements}

\author{Christos~A.~Athanasiadis}
\address{Department of Mathematics (Division of Algebra-Geometry)\\
University of Athens\\
Panepistimioupolis\\
15784 Athens, Greece}
\email{caath@math.uoa.gr}

\date{October 16, 2006}
\begin{abstract}
Given a nonnegative integer $m$ and a finite collection $\aA$ of
linear forms on $\QQ^d$, the arrangement of affine hyperplanes in
$\QQ^d$ defined by the equations $\alpha(x) = k$ for $\alpha \in \aA$
and integers $k \in [-m, m]$ is denoted by $\aA^m$. It is proved that
the coefficients of the characteristic polynomial of $\aA^m$ are
quasi-polynomials in $m$ and that they satisfy a simple combinatorial
reciprocity law.
\end{abstract}

\maketitle

\section{Introduction}
\label{intro}

Let $V$ be a $d$-dimensional vector space over the field $\QQ$ of rational
numbers and $\aA$ be a finite collection of linear forms on $V$ which
spans the dual vector space $V^*$. We denote by $\aA^m$ the essential
arrangement of affine hyperplanes in $V$ defined by the equations
$\alpha(x) = k$ for $\alpha \in \aA$ and integers $k \in [-m, m]$ (we
refer to \cite{OT, Sta3} for background on hyperplane arrangements).
Thus $\aA^0$ consists of the linear hyperplanes which are the kernels of
the forms in $\aA$ and $\aA^m$ is a deformation of $\aA^0$, in the sense
of \cite{Ath1, PS}.

The characteristic polynomial \cite[Section 2.3]{OT} \cite[Section
1.3]{Sta3} of $\aA^m$, denoted $\chi_\aA (q, m)$, is a fundamental
combinatorial and topological invariant which can be expressed as
\begin{equation}
\chi_\aA (q, m) \, = \, \sum_{i=0}^d \ c_i (m) \, q^i.
\label{eq:qm}
\end{equation}
We will be concerned with the behavior of $\chi_\aA (q, m)$ as a
function of $m$. Let $\NN := \{0, 1,\dots\}$ and recall that a
function $f: \NN \to \RR$ is called a \emph{quasi-polynomial} with
period $N$ if there exist polynomials $f_1, f_2,\dots,f_N: \NN \to
\RR$ such that $f(m) = f_i (m)$ for all $m \in \NN$ with $m \equiv
i \ ({\rm mod} \, N)$. The degree of $f$ is the maximum of the
degrees of the $f_i$. Our main result is the following theorem.

\begin{theorem}
Under the previous assumptions on $\aA$, the coefficient $c_i (m)$
of $q^i$ in $\chi_\aA (q, m)$ is a quasi-polynomial in $m$ of degree
at most $d-i$. Moreover, the degree of $c_0 (m)$ is equal to $d$ and
\begin{equation}
\chi_\aA (q, -m) \, = \, (-1)^d \chi_\aA (-q, m-1).
\label{eq:rec}
\end{equation}
\label{thm0}
\end{theorem}

In particular we have $\chi_\aA (q, -1) = (-1)^d \chi_\aA (-q)$,
where $\chi_\aA (q)$ is the characteristic polynomial of $\aA^0$.
Let $\aA_\RR^m$ denote the arrangement of affine hyperplanes in the
real $d$-dimensional vector space $V_\RR = V \otimes_\QQ \RR$ defined
by the same equations defining the hyperplanes of $\aA^m$. Let $r_\aA
(m) = (-1)^d \chi_\aA (-1, m)$ and $b_\aA (m) = (-1)^d \chi_\aA (1,
m)$ so that, for $m \in \NN$, $r_\aA (m)$ and $b_\aA (m)$ count the
number of regions and bounded regions, respectively, into which $V_\RR$
is dissected by the hyperplanes of $\aA_\RR^m$ \cite[Section 2.2]{Sta3}
\cite{Za}.

\begin{corollary}
Under the previous assumptions on $\aA$, the function $r_\aA (m)$
is a quasi-polynomial in $m$ of degree $d$ and, for all positive
integers $m$, $(-1)^d r_\aA (-m)$ is equal to the number $b_\aA (m-1)$
of bounded regions of $\aA_\RR^{m-1}$.
\label{cor0}
\end{corollary}

Theorem \ref{thm0} and its corollary belong to a family of results
demonstrating some kind of combinatorial reciprocity law; see \cite{Sta1}
for a systematic treatment of such phenomena.
Not surprisingly, the proof given in Section \ref{proof} is a simple
application of the main results of Ehrhart theory \cite[Section 4.6]{Sta2}.
More specifically, equation (\ref{eq:rec}) will follow from the reciprocity
theorem \cite[Theorem 4.6.26]{Sta2} for the Ehrhart quasi-polynomial of
a rational polytope. An expression for the coefficient of the leading
term $m^d$ of either $c_0 (m)$ or $r_\aA (m)$ is also derived in that
section. Some examples, including the motivating example in which
$\aA_\RR^0$ is the arrangement of reflecting hyperplanes of a Weyl group,
and remarks are discussed in Section \ref{remarks}. In the remainder of
this section we give some background on characteristic and Ehrhart
(quasi-)polynomials needed in Section \ref{proof}. We will denote by $\#
S$ or $|S|$ the cardinality of a finite set $S$.

\bigskip
\noindent
\textit{Arrangements of hyperplanes.}
Let $V$ be a $d$-dimensional vector space over a field $\KK$. An
\emph{arrangement of hyperplanes} in $V$ is a finite collection $\hH$ of
affine subspaces of $V$ of codimension one (we will allow this collection to
be a multiset). The \emph{intersection poset} of $\hH$ is the set $L_\hH =
\{ \cap \, \fF: \fF \subseteq \hH\}$ of all intersections of subcollections
of $\hH$, partially ordered by reverse inclusion. It has a unique minimal
element $\hat{0} = V$, corresponding to the subcollection $\fF = \emptyset$.
The \emph{characteristic polynomial} of $\hH$ is defined by
\begin{equation}
\chi_\hH (q) \, = \sum_{x \in L_\hH} \mu (x) \, q^{\dim x}
\label{eq:char}
\end{equation}
where $\mu$ stands for the M\"obius function on $L_\hH$ defined by
\[ \mu (x) \, = \, \begin{cases}
1, & \text{if \ $x = \hat{0}$} \\
- \sum_{y < x} \mu (y), & \text{otherwise.} \end{cases} \]
Equivalently \cite[Lemma 2.55]{OT} we have
\begin{equation}
\chi_\hH (q) \, = \sum_{\gG \subseteq \hH} (-1)^{\# \gG} \,
q^{\dim (\cap \, \gG)}
\label{eq:char2}
\end{equation}
where the sum is over all $\gG \subseteq \hH$ with $\cap \, \gG \neq
\emptyset$.

In the case $\KK = \RR$, the connected components of the space obtained
from $V$ by removing the hyperplanes of $\hH$ are called \emph{regions}
of $\hH$. A region is \emph{bounded} if it is a bounded subset of $V$
with respect to a usual Euclidean metric.

\bigskip
\noindent
\textit{Ehrhart quasi-polynomials.}
A convex polytope $P \subseteq \RR^n$ is said to be a \emph{rational} or
\emph{integral} polytope if all its vertices have rational or integral
coordinates, respectively. If $P$ is rational and $P^\circ$ is its
relative interior then the functions defined for nonnegative integers
$m$ by the formulas
\begin{equation}
\begin{tabular}{l}
$i (P, m) \, = \, \# \, (m P \cap \ZZ^n)$ \\
$\bar{i} (P, m) \, = \, \# \, (m P^\circ \cap \ZZ^n)$
\end{tabular}
\label{eq:ehr}
\end{equation}
are quasi-polynomials in $m$ of degree $d = \dim (P)$, related by the
Ehrhart reciprocity theorem \cite[Theorem 4.6.26]{Sta2}
\begin{equation}
i (P, -m) \, = \, (-1)^d \ \bar{i} (P, m).
\label{eq:ehrrec}
\end{equation}
The function $i (P, m)$ is called the \emph{Ehrhart quasi-polynomial} of
$P$. The coefficient of the leading term $m^d$ in either $i (P, m)$ or
$\bar{i} (P, m)$ is a constant equal to the normalized $d$-dimensional
volume of $P$ (meaning the $d$-dimensional volume of $P$ normalized with
respect to the affine lattice $V_P \cap \ZZ^n$, where $V_P$ is the affine
span of $P$ in $\RR^n$). If $P$ is an integral polytope then $i (P, m)$
is a polynomial in $m$ of degree $d$, called the \emph{Ehrhart
polynomial} of $P$.

\section{Proof of Theorem \ref{thm0}}
\label{proof}

In this section we prove Theorem \ref{thm0} and Corollary \ref{cor0}
and derive a formula for the coefficient of the leading term $m^d$ of
$r_\aA (m)$. In what follows $\aA$ is as in the beginning of Section
\ref{intro}. We use the notation $[a, b] = \{x \in \RR: a \le x \le b\}$
and $[a, b]_\ZZ = [a, b] \cap \ZZ$ for $a, b \in \ZZ$ with $a \le b$.

\vspace{0.1 in}
\noindent
\emph{Proof of Theorem \ref{thm0} and Corollary \ref{cor0}.} Using
formula (\ref{eq:char2}) we get
\begin{equation}
\chi_{\aA^m} (q) \, = \sum_{\gG \subseteq \aA^m} (-1)^{\# \gG} \,
q^{\dim (\cap \, \gG)}
\label{eq:proof1}
\end{equation}
where the sum is over all $\gG \subseteq \aA^m$ with $\cap \, \gG
\neq \emptyset$. Clearly for this to happen $\gG$ must contain at most
one hyperplane of the form $\alpha (x) = k$ for each $\alpha \in \aA$.
In other words we must have $\gG = \fF_b$ for some $\fF \subseteq \aA$
and map $b: \fF \to [-m,m]_\ZZ$ sending $\alpha$ to $b_\alpha$, where
$\fF_b$ consists of the hyperplanes $\alpha (x) = b_\alpha$ for
$\alpha \in \fF$. Let us denote by $\dim \fF$ the dimension of the
linear span of $\fF$ in $V^*$ and observe that $\dim (\cap \, \fF_b)
= d - \dim \fF$ whenever $\cap \, \fF_b$ is nonempty. From the
previous observations and (\ref{eq:proof1}) we get

\begin{eqnarray*}
\chi_\aA (q, m) & = & \sum_{\fF \subseteq \aA} \ \sum_{\substack{b:
\fF \to [-m,m]_\ZZ \\ \cap \, \fF_b \neq \emptyset}}
(-1)^{\# \fF_b} \, q^{\dim (\cap \, \fF_b)} \\ \\
& = & \sum_{\fF \subseteq \aA} (-1)^{\# \fF} \,
q^{d-\dim \fF} \ \# \{b: \fF \to [-m,m]_\ZZ, \, \cap \, \fF_b \neq
\emptyset \}.
\end{eqnarray*}
Let us write $\fF = \{\alpha_1, \alpha_2,\dots,\alpha_n\}$ and $b_i
= b_{\alpha_i}$, so that $b$ can be identified
with a column vector in $\QQ^n$. Then $\cap \, \fF_b$ is nonempty
if and only if the linear system $\alpha_i (x) = b_i$, $1 \le i \le
n$, has a solution in $\QQ^d$ or, equivalently, if and only if $b$
lies in the image $\im T_\fF$ of the linear transformation $T_\fF:
\QQ^d \to \QQ^n$ mapping $x \in \QQ^d$ to the column vector in
$\QQ^n$ with coordinates $\alpha_1 (x)$, $\alpha_2
(x),\dots,\alpha_n (x)$. It follows that

\begin{eqnarray*}
\# \{b: \fF \to [-m,m]_\ZZ, \, \cap \, \fF_b \neq \emptyset \} & =
& \# \ \im T_\fF \cap ([-m,m]_\ZZ)^n \\
& = & \# \ \im T_\fF \cap [-m,m]^n \cap \ZZ^n \\
& = & \# \ (m \, (\im T_\fF \cap [-1,1]^n) \cap \ZZ^n) \\
& = & \# \ (m P_\fF \cap \ZZ^n) \\
& = & i (P_\fF, m)
\end{eqnarray*}

\noindent
where $P_\fF = (\im T_\fF \otimes_\QQ \RR) \cap [-1,1]^n$, and hence
that
\begin{equation}
\chi_\aA (q, m) \, = \sum_{\fF \subseteq \aA} (-1)^{\# \fF} \, q^{d-\dim
\fF} \ i (P_\fF, m).
\label{eq:proof2}
\end{equation}
Equivalently we have
\begin{equation}
c_i (m) \ = \sum_{\substack{\fF \subseteq \aA \\ \dim \fF = d-i}}
(-1)^{\# \fF} \, i (P_\fF, m)
\label{eq:proof3}
\end{equation}
for $0 \le i \le d$, where the $c_i (m)$ are as in (\ref{eq:qm}).
Clearly $P_\fF$ is a rational convex polytope of dimension $\dim (\im
T_\fF) = \dim \fF$ and hence $i (P_\fF, m)$ is a quasi-polynomial in
$m$ of degree $\dim \fF$. It follows from (\ref{eq:proof3}) that
$c_i (m)$ is a quasi-polynomial in $m$ of degree at most $d-i$ and
that $r_\aA (m) = \sum_{i=0}^d \, (-1)^{d-i} c_i (m)$ is a
quasi-polynomial in $m$ of degree at most $d$. Moreover we have $r_\aA
(m) \ge (2m+2)^d$ for $m \ge 0$ since $\aA$ contains $d$ linearly
independent forms and the corresponding hyperplanes of $\aA^m_\RR$
dissect $V_\RR$ into $(2m+2)^d$ regions. It follows that the degree
of $r_\aA (m)$ is no less than $d$, which implies that the degrees
of $r_\aA (m)$ and $c_0 (m)$ are, in fact, equal to $d$.

It remains to prove the reciprocity relation (\ref{eq:rec}). For
$\fF \subseteq \aA$ with $\# \fF = n$ let $W_\fF$ be the real linear
subspace $\im T_\fF \otimes_\QQ \RR$ of $\RR^n$, so that $P_\fF =
W_\fF \cap [-1,1]^n$. We have

\begin{eqnarray*}
m P^\circ_\fF \cap \ZZ^n & = & (W_\fF \cap [-m,m]^n)^\circ \cap
\ZZ^n \\
& = & W_\fF \cap [-(m-1), m-1]^n \cap \ZZ^n \\
& = &  (m-1) P_\fF \cap \ZZ^n
\end{eqnarray*}
and hence $\bar{i} (P_\fF, m) = i (P_\fF, m-1)$. The Ehrhart
reciprocity theorem (\ref{eq:ehrrec}) implies that
\begin{equation}
i (P_\fF, -m) \, = \, (-1)^{\dim \fF} \, i (P_\fF, m-1).
\label{eq:proof4}
\end{equation}
Equation (\ref{eq:rec}) follows from (\ref{eq:proof2}) and
(\ref{eq:proof4}).
\qed

\bigskip
The following corollary is an immediate consequence of the case
$i=0$ of (\ref{eq:proof3}), the equation $r_\aA (m) = \sum_{i=0}^d
\, (-1)^{d-i} c_i (m)$ and the fact that the degree of $c_i (m)$
is less than $d$ for $1 \le i \le d$.
\begin{corollary}
The coefficient of the leading term $m^d$ in $r_\aA (m)$ is equal
to the expression
\[ \sum_{\substack{\fF \subseteq \aA \\ \dim \fF = d}}
(-1)^{\# \fF - d} \, \vol_d (P_\fF), \]
where $P_\fF$ is as in the proof of Theorem \ref{thm0} and $\vol_d
(P_\fF)$ is the normalized $d$-dimensional volume of $P_\fF$. \qed
\label{cor1}
\end{corollary}

\section{Examples and remarks}
\label{remarks}

In this section we list a few examples, questions and remarks.
\begin{example} {\rm
If $V = \QQ$ and $\aA$ consists of two forms $\alpha_1, \alpha_2: V \to
\QQ$ with $\alpha_1 (x) = x$ and $\alpha_2 (x) = 2x$ for $x \in V$ then
$\aA^m$ consists of the affine hyperplanes (points) in $V$ defined by the
equations $x = k$ and $x = k/2$ for $k \in [-m, m]_\ZZ$. One can check
that
\[ \chi_\aA (q, m) \ = \
\begin{cases} q-3m-1, & \text{if $m$ is even} \\
q-3m-2, & \text{if $m$ is odd} \end{cases} \]
and that (\ref{eq:rec}) holds. Moreover we have
\[ r_\aA (m) \ = \
\begin{cases} 3m+2, & \text{if $m$ is even} \\
3m+3, & \text{if $m$ is odd}. \end{cases} \]
Note that $\vol_d (P_\fF)$ takes the values 2, 2 and 1 when
$\fF = \{\alpha_1\},
\{\alpha_2\}$ and $\{\alpha_1, \alpha_2\}$, respectively.
}
\label{ex1}
\end{example}
\begin{example} {\rm
If $V = \QQ^d$ and $\aA$ consists of the coordinate functions
$\alpha_i (x) = x_i$ for $1 \le i \le d$ then $\aA^m$ consists of
the affine hyperplanes in $V$ given by the equations $x_i = k$ with
$1 \le i \le d$, $k \in [-m, m]_\ZZ$ and $\chi_\aA (q, m) = (q-2m-1)^d$,
which is a polynomial in $q$ and $m$ satisfying (\ref{eq:rec}).
}
\label{ex2}
\end{example}
\begin{example} {\rm
Let $\Phi$ be a finite, irreducible, crystallographic root system
spanning the Euclidean space $\RR^d$, endowed with the standard inner
product $( \ , \, )$ (we refer to \cite{BB, Bou, Hu} for background on
root systems). Fix a positive system $\Phi^+$ and let $Q_\Phi$ and $W$
be the coroot lattice and Weyl group, respectively, corresponding to
$\Phi$. Let also $\aA_\Phi^m$ denote the $m$th \emph{generalized
Catalan arrangement}
associated to $\Phi$ \cite{Ath1, Ath2, PS}, consisting of the affine
hyperplanes in $\RR^d$ defined by the equations $(\alpha, x) = k$ for
$\alpha \in \Phi^+$ and $k \in [-m, m]_\ZZ$ (so that $\aA_\Phi^0$ is
the real reflection arrangement associated to $\Phi$). If $V$ is the
$\QQ$-span
of $Q_\Phi$ then there exists a finite collection $\aA$ of linear forms
on $V$ (one for each root in $\Phi^+$) such that, in the notation of
Section \ref{intro}, $\aA_\RR^m$ coincides with $\aA_\Phi^m$. In
\cite[Theorem 1.2]{Ath2} a uniform proof was given of the formula
\begin{equation}
\chi_\aA (q, m) \, = \, \prod_{i=1}^d \, (q-mh-e_i)
\label{eq:mcat}
\end{equation}
for the characteristic polynomial of $\aA_\Phi^m$, where $e_1,
e_2,\dots,e_d$ are the exponents and $h$ is the Coxeter number of
$\Phi$. Thus the reciprocity law (\ref{eq:rec}) in this case is
equivalent to the well-known fact \cite[Section V.6.2]{Bou}
\cite[Lemma 3.16]{Hu} that the numbers $h - e_i$ are a permutation
of the $e_i$. As was already deduced in \cite[Corollary 1.3]{Ath2},
it follows from (\ref{eq:mcat}) that
\[ r_\aA (m) \, = \, \prod_{i=1}^d \, (mh+e_i+1) \]
and
\[ b_\aA (m) \, = \, \prod_{i=1}^d \, (mh+e_i-1) \]
are polynomials in $m$ of degree $d$ (a fact which was the main
motivation behind this paper). Setting $N (\Phi, m) =
\frac{1} {|W|} \, r_\aA (m)$ and $N^+ (\Phi, m) = \frac{1} {|W|}
\, b_\aA (m)$, as in \cite{AT, FR}, our Corollary \ref{cor0} implies
that
\begin{equation}
(-1)^d \, N(\Phi, -m) \, = \, N^+ (\Phi, m-1).
\label{eq:N}
\end{equation}
It was suggested in \cite[Remark 12.5]{FR} that this equality, first
observed in \cite[(2.12)]{FR}, may be an instance of Ehrhart
reciprocity. This was confirmed in \cite[Section 7]{AT} using an
approach which is different from the one followed in this paper.
Finally we note that Corollary \ref{cor1} specializes to the
curious identity
\begin{equation}
h^d \, = \, \sum_F \, (-1)^{\# F - d} \, \vol_d (P_F)
\label{eq:curious}
\end{equation}
where in the sum on the right hand-side $F$ runs through all
subsets $\{\alpha_1, \alpha_2,\dots,\alpha_n\}$ of $\Phi^+$
spanning $\RR^d$, $P_F$ is the intersection of the cube $[-1,
1]^n$ with the image of the linear transformation $T_F: \RR^d \to
\RR^n$ mapping $x \in \RR^d$ to the column vector in $\RR^n$ with
coordinates $(\alpha_1, x)$, $(\alpha_2, x),\dots,(\alpha_n, x)$
and $\vol_d (P_F)$ is the normalized $d$-dimensional volume of
$P_F$. If $\Phi$ has type $A_d$ in the Cartan-Killing
classification then (\ref{eq:curious}) translates to the equation
\begin{equation}
(d+1)^d \, = \, \sum_G \, (-1)^{e(G) - d} \, \vol_d (Q_G)
\label{eq:curiousA}
\end{equation}
where in the sum on the right hand-side $G$ runs through all
connected simple graphs on the vertex set $\{1, 2,\dots,d+1\}$,
$e(G)$ is the number of edges of $G$ and $Q_G$ is the
$d$-dimensional polytope in $\RR^d$ defined in the following way.
Let $\tau$ be a spanning tree of $G$ with edges labelled in a one
to one fashion with the variables $x_1, x_2,\dots,x_d$. For any
edge $e$ of $G$ which is not an edge of $T$ let $R_e$ be the
region of $\RR^d$ defined by the inequalities $-1 \le x_{i_1} +
x_{i_2} + \cdots + x_{i_k} \le 1$, where $x_{i_1},
x_{i_2},\dots,x_{i_k}$ are the labels of the edges (other than
$e$) of the fundamental cycle of the graph obtained from $T$ by
adding the edge $e$. The polytope $Q_G$ is the intersection of the
cube $[-1, 1]^d$ and the regions $R_e$. } \label{ex3}
\end{example}
\begin{remark} {\rm
It is well-known \cite[Corollary 3.5]{Sta3} that the coefficients
of the characteristic polynomial of a hyperplane arrangement strictly
alternate in sign. As a consequence, in the notation of (\ref{eq:qm}),
we have $(-1)^{d-i} c_i (m) > 0$ for all $m \in \NN$ and $0 \le i \le
d$. We do not know of an example of a collection $\aA$ of forms for a
which a negative number appears among the coefficients of the
quasi-polynomials $(-1)^{d-i} c_i (m)$.
}
\label{rem0}
\end{remark}
\begin{remark} {\rm
If the matrix defined by the forms in $\aA$ with respect to some
basis of $V$ is integral and totally unimodular, meaning that all
its minors are $-1, 0$ or 1, then the polytopes $P_\fF$ in the
proof of Theorem \ref{thm0} are integral and, as a consequence,
the functions $c_i (m)$ and $r_\aA (m)$ are polynomials in $m$.
This assumption on $\aA$ is satisfied in the case of graphical
arrangements, that is when $\aA$ consists of the forms $x_i - x_j$
on $\QQ^r$, where $1 \le i < j \le r$, corresponding to the edges
$\{i, j\}$ of a simple graph $G$ on the vertex set $\{1,
2,\dots,r\}$. The degree of the polynomial $r_G (m) := r_\aA (m)$
is equal to the dimension of the linear span of $\aA$, in other
words to the rank of the cycle matroid of $G$.
}
\label{rem1}
\end{remark}
\begin{remark} {\rm
Let $\aA$ and $\hH$ be finite collections of linear forms on a
$d$-dimensional $\QQ$-vector space $V$ spanning $V^*$. Using
the notation of Section \ref{intro}, let $\hH_m$ denote the union
of $\aA_\RR^m$ with the linear arrangement $\hH_\RR^0$. It follows
from Theorem \ref{thm0}, the
Deletion-Restriction theorem \cite[Theorem 2.56]{OT} and induction
on the cardinality of $\hH$ that the function $r (\hH_m)$ is a
quasi-polynomial in $m$ of degree $d$.
Given a region $R$ of $\hH_\RR^0$, let $r_R (m)$ denote the
number of regions of $\hH_m$ which are contained in $R$, so that
\[ r (\hH_m) \, = \, \sum_R \, r_R (m)\]
where $R$ runs through the set of all regions of $\hH_\RR^0$. Is the
function $r_R (m)$ always a quasi-polynomial in $m$?
}
\label{rem2}
\end{remark}

\end{document}